\documentclass[12pt]{article}%
\usepackage{amsfonts}
\usepackage{amsmath,amssymb,amsthm,enumerate,epsfig,graphicx,natbib, ulem}
\usepackage{ifthen,latexsym,syntonly}
\usepackage{amsmath}
\usepackage{amssymb}
\usepackage{graphicx}
\usepackage{verbatim}
\usepackage{morefloats}
\RequirePackage{natbib}
\usepackage{natbib}
\usepackage{multibib}
\newcites{app}{References}
\usepackage{geometry}
\usepackage{soul,color}%
\usepackage{cprotect}
\usepackage{paralist}
\usepackage{float}
\usepackage{booktabs}
\usepackage{bbm}

\newtheorem{thm}{Theorem}

\newtheorem{prop}{Proposition}

\newtheorem{rmk}{Remark}

\newtheorem{assump}{Assumption}
\newenvironment{myassump}[2][]
{\renewcommand\theassump{#2}
	\begin{assump}[#1]}
	{\end{assump}}

 \DeclareMathOperator{\cov}{Cov}
\DeclareMathOperator{\tr}{tr}

\def\ul{\underline}

\def\fZ {{\mathbf Z}}
\def\fY {{\mathbf Y}}

\def\fS {{\mathbf S}}

\def\fr {{\mathbf r}}

\def\fA {{\mathbf A}}
\def\fB {{\mathbf B}}

\def\fI {{\mathbf I}}

\def\fQ {{\mathbf Q}}
\def\fSigma {{\boldsymbol{\Sigma}}}

\def\ff {{\mathbf f}}

\newcommand{\ra}[1]{\renewcommand{\arraystretch}{#1}}

\makeatletter
\renewcommand{\tablename}{{\bf Table}}
\renewcommand{\fnum@table}[1]{\normalfont\textbf{\tablename~\thetable}\\}

\renewcommand{\theenumi}{\Alph{enumi}}

 \makeatletter
 \renewcommand{\p@enumii}{\theenumi.}
 \makeatother

\makeatother
\usepackage{indentfirst} 

\makeatletter
\renewcommand{\figurename}{{\bf Fig.}}
\renewcommand{\fnum@figure}[1]{\normalfont\textbf{\figurename~\thefigure.}}
\makeatother

\makeatletter
\newif\ifrhsapp
\rhsappfalse

\newcommand{\@seccntformat@section}[1]{%
  \ifrhsapp
  Appendix
  \else
  \fi
  \csname the#1\endcsname.\quad
}

\newcommand*{\@seccntformat@subsection}[1]{%
  \csname the#1\endcsname.\quad
}

\newcommand*{\@seccntformat@subsubsection}[1]{%
  \csname the#1\endcsname.\quad
}

\let\@@seccntformat\@seccntformat
\renewcommand*{\@seccntformat}[1]{%
  \expandafter\ifx\csname @seccntformat@#1\endcsname\relax
    \expandafter\@@seccntformat
  \else
    \expandafter
      \csname @seccntformat@#1\expandafter\endcsname
  \fi
    {#1}%
}

\renewcommand{\subsection}{\@startsection
  {subsection}{2}{0mm}{-3.25ex \@plus -1ex \@minus -.2ex}{1.5ex \@plus.2ex}{\normalfont\large\itshape}}

\renewcommand{\subsubsection}{\@startsection
  {subsubsection}{2}{0mm}{-3.25ex \@plus -1ex \@minus -.2ex}{1.5ex \@plus.2ex}{\normalfont\itshape}}

\def\appendix{\par
  \setcounter{section}{0}%
  \setcounter{subsection}{0}%
  \rhsapptrue
  \renewcommand\thesection{\Alph{section}}%
}
\makeatother

\begin{document}

\title{Testing High-dimensional Covariance Matrices under the Elliptical Distribution and Beyond
\thanks{
Correspondence to: Department of Information Systems, Business Statistics and Operations Management, Hong Kong University of Science and Technology, Clear Water Bay, Kowloon, Hong Kong. Tel:(+852) 2358 7750. E-mail address: xhzheng@ust.hk
}}
\author{
\sc  Xinxin Yang$^1$, Xinghua Zheng$^{2,*}$, Jiaqi Chen$^3$\\
       {\footnotesize $ ^{1}$
       School of Statistics and Mathematics, Central University of Finance and Economics}\\
       {\footnotesize $ ^{2}$
       Department of ISOM, Hong Kong University of Science and Technology,}\\
        {\footnotesize $ ^{3}$
       Department of Mathematics, Harbin Institute of Technology}\\
}
\date{}
\maketitle

\begin{abstract}
We develop tests for high-dimensional covariance matrices under a generalized elliptical model. Our tests are based on a central limit theorem (CLT) for  linear spectral statistics of the sample covariance matrix based on self-normalized observations. For testing sphericity, our tests neither assume specific parametric distributions nor involve the kurtosis of data. More generally, we can test against any non-negative definite matrix that can even be not invertible. As an interesting application, we illustrate in empirical studies that our tests can be used to test uncorrelatedness among idiosyncratic returns.\vskip 0.8cm
\noindent {\bf Keywords:}  Covariance matrix, high-dimension, elliptical model, linear spectral statistics, central limit theorem.


\end{abstract}

\newpage


\section{Introduction}
\label{sec:intro}

\subsection{Tests for high-dimensional covariance matrices}\label{test_hd}
Testing covariance matrices is of fundamental importance in multivariate analysis.
 There has been a long history of study on testing (i) the covariance matrix $\boldsymbol{\Sigma}$ is equal to a given matrix, or (ii) the covariance matrix~$\boldsymbol{\Sigma}$ is proportional to a given matrix. Specifically, for a given non-negative definite matrix $\boldsymbol{\Sigma}_0$, one aims to {test}
\begin{align}\label{test:identity}
	&H_{0}:~\boldsymbol{\Sigma} = \boldsymbol{\Sigma}_0\quad vs. \quad H_{a}:~\boldsymbol{\Sigma} \neq\boldsymbol{\Sigma}_0,
\end{align}
or
\begin{align}
	&H_{0}:~\boldsymbol{\Sigma} \propto \boldsymbol{\Sigma}_0\quad vs. \quad H_{a}:~\boldsymbol{\Sigma} \not\propto\boldsymbol{\Sigma}_0,\label{test:sphericity}
	\end{align}
where ``$\propto$" stands for ``proportional to". When $\fSigma_0=\fI$, test \eqref{test:identity} is referred to as the identity test and \eqref{test:sphericity} as the sphericity  test. If $\fSigma_0$ is invertible, then testing \eqref{test:identity} or \eqref{test:sphericity} can be reduced to the identity or sphericity test, by multiplying the observations with $\fSigma_0^{-1/2}$.

In the classical setting where the dimension $p$ is fixed and the sample size $n$ goes to infinity, the sample covariance matrix is a consistent estimator, and further inference can be made based on the associated {\it central limit theory} (CLT). Examples include the likelihood ratio tests (see, e.g., \cite{Muirhead82}, Sections~8.3 and 8.4), and the locally most powerful invariant tests (\cite{John1971}, \cite{Nagao73}).

In the high-dimensional setting, because the sample covariance matrix is inconsistent, conventional tests may not apply. For the identity and sphericity tests, new methods have been developed, first under the multivariate normal distribution, then for more generally distributed data:
\begin{itemize}
	\item Multivariate normally distributed data. When $p/n\rightarrow y\in(0, \infty)$, \cite{LedoitW02} show that John's test of sphericity is still consistent and propose a modified Nagao's identity test. \cite{Srivastava05} introduces a new test of sphericity under a more general condition that $n=O(p^{\delta})$ for some $\delta\in(0, 1]$. \cite{Birke05} show that the asymptotic null distributions of John's and the modified Nagao's test statistics in \cite{LedoitW02} are still valid when $p/n\rightarrow\infty$. Relaxing the normality assumption but still assuming the kurtosis equals $3$, \cite{BaiJ09} develop a corrected likelihood ratio  test of identity when $p/n\rightarrow y\in(0, 1)$. {For testing sphericity, \cite{JiangYAos13} derive the asymptotic distribution of  the likelihood ratio  test statistic under the multivariate normal distribution with $p/n\rightarrow y\in(0, 1]$}.
\item More generally distributed data.  \cite{ChenZ10}  generalize the results in \cite{LedoitW02} without assuming normality nor an explicit relationship between $p$ and $n$. By relaxing the kurtosis assumption, \cite{WCY13} extend the corrected likelihood ratio test in \cite{BaiJ09} and the modified Nagao's test in \cite{LedoitW02}. Along this line, \cite{WangY13} propose two tests by correcting the likelihood ratio test and John's test.
\end{itemize}	

\subsection{The elliptical distribution and its applications}
The elliptically distributed data can be expressed as
$$
\fY=\omega\fSigma^{1/2}\fZ,
$$
where $\omega$ is a positive random scalar, $\fZ$ is a $p$-dimensional random vector from $N(\mathbf{0}, \fI)$, and further~$\omega$ and $\fZ$ are independent of each other. It is a natural generalization of the multivariate normal distribution, and contains many widely used
distributions as special cases including the multivariate $t$-distribution, the symmetric multivariate Laplace distribution and the symmetric multivariate stable distribution.
See \cite{FangK90} for further details.

One of our motivations of this study arises from the wide applicability of the elliptical distribution. The ``mixture coefficient'' $\omega$ can feature heteroskedasticity that are widely present in real data. Furthermore, just like the ARCH (\cite{Engle82}) and GARCH (\cite{Bollerslev86}) models, the marginal distribution of $\fY$ is a mixture of normal hence heavy-tailed. Therefore, the elliptical distribution can feature both heteroskedasticity and heavy-tailedness. In finance,  stock returns have been extensively documented to exhibit such two features, dating back at least to \cite{fama1965behavior} and \cite{mandelbrot1967variation}.
{Accommodating heteroskedasticity and heavy-tailedness makes} the elliptical distribution a more admissible candidate for stock-return models than the Gaussian distribution; see, e.g., \cite{owen83} and \cite{bingham2002semi}.
\cite{mcneil05} state that ``elliptical distributions ... provided far superior models to the
multivariate normal for daily and weekly US stock-return data" and that ``multivariate return data for groups of returns of similar type often look roughly elliptical."

\subsection{Performance of  existing tests under the elliptical model} \label{performance_of_existing_tests}
Given the wide applicability of the elliptical distribution, it is important to check whether existing tests for covariance matrices are applicable to the elliptical distribution. Both numerical and theoretical analysis give a negative answer.

We start with a simple numerical study to investigate the empirical sizes.  Consider observations
$\fY_i=\omega_i \fZ_i, i=1, \cdots, n,$ {where}
\begin{enumerate}[(i)]
	\item $\omega_i$'s are  absolute values of \hbox{i.i.d.} standard normal random variables,
	\item $\fZ_i$'s are \hbox{i.i.d.} $p$-dimensional standard multivariate normal random vectors, and
	\item  $\omega_i$'s and $\fZ_i$'s are independent of each other.
\end{enumerate}
Under such a setting, $\fY_i$'s are \hbox{i.i.d.} random vectors with mean $\mathbf{0}$ and covariance matrix $\fI$. We will test both  $H_0: \boldsymbol{\Sigma}=\fI$ and $H_0: \boldsymbol{\Sigma}\propto\fI$.

To test $H_0: \boldsymbol{\Sigma}=\fI$, we use the tests {in} \cite{LedoitW02} (LW$_1$ test), \cite{BaiJ09} (BJYZ test),  \cite{ChenZ10} (CZZ$_1$ test) and \cite{WCY13} (WYMC-LR and WYMC-LW tests). For testing $H_0: \boldsymbol{\Sigma}\propto\fI$, we apply the tests proposed in \cite{LedoitW02} (LW$_2$ test), \cite{Srivastava05} (S test),  \cite{ChenZ10} (CZZ$_2$ test) and \cite{WangY13} (WY-LR and WY-JHN tests).
Table \ref{counterexample_sphere} reports the empirical sizes for these tests at~$5\%$ significance~level.

\begin{center}
+++  Insert Table \ref{counterexample_sphere} Here +++
\end{center}

We observe from Table \ref{counterexample_sphere} that the empirical sizes of all these tests are far higher than the nominal level of $5\%$,
suggesting that they are inconsistent  under the elliptical distribution.

It is worth discussing why these tests fail. Firstly, the reason is not that these tests do not apply to the high-dimensional setting. In fact, when $\omega_i \equiv 1$, for the same pairs of dimensions and sample sizes, all these tests yield  sizes close to the nominal level of 5\%; see Table \ref{compare_size} for details. Secondly, it is not due to heavy-tailedness either. The marginal distribution of~$\fY$, although is heavier than normal, still has exponentially decaying tails with a finite moment generating function.

The real reason that these tests fail lies in the presence of $\{\omega_i\}$. Denote $\fS_n=n^{-1}\sum_{i=1}^n\fY_i\fY_i^{\top}=n^{-1}\sum_{i=1}^n\omega_i^2 \fZ_i\fZ_i^{\top}$. The celebrated Mar\v{c}enko-Pastur theorem states that, \emph{when $\omega_i$'s are constant}, the {\it empirical spectral distribution} (ESD) of $\fS_n$ converges to the Mar\v{c}enko-Pastur law. This convergence turns out to be crucial for all the aforementioned tests as they all involve certain moments of the limiting ESD (LSD). When $\omega_i$'s are \emph{not constant}, Theorem 1 of \cite{ZhengL11} implies that the ESD of $\fS_n$ will \emph{not} converge to the Mar\v{c}enko-Pastur law. Consequently, the asymptotic null distributions of the aforementioned test statistics change, and the tests no longer apply.

\subsection{Our model and aim of this study}\label{Ourmodelandgoal}
The previous section shows that  existing tests do not apply when observations are heteroskedastic, a feature that is commonly encountered in finance, economics and many other fields.


In this paper, we study tests for high-dimensional covariance matrices when  data may exhibit heteroskedasticity. Specifically, we consider the following model.  Denote by $\fY_i$, $i=1,\cdots, n$, the observations, which can be written as
\begin{equation}\label{observe}
\fY_i=\omega_i\fSigma^{1/2}\fZ_i,
\end{equation}
where
\begin{enumerate}[(i)]
	\item $\omega_i$'s are positive random scalars reflecting heteroskedasticity,
	\item $\fSigma\in\mathbb{R}^{p\times p}$ is a non-negative definite matrix,
	\item {$\fZ:=\big(\fZ_1, \ldots, \fZ_n\big)=(Z_{ij})_{p\times n}$ consists of \hbox{i.i.d.} standardized random variables,}
	\item { $\omega_i$'s can depend on each other and on $\{\fZ_i: ~i=1,\cdots, n\}$ in an \emph{arbitrary} way, and}
	\item $\omega_i$'s do \emph{not} need to be stationary.
\end{enumerate}

Model \eqref{observe} incorporates the elliptical distribution as a special case. This general model further possesses several important advantages:
\begin{itemize}
	\item It can be considered as a multivariate extension of the ARCH/GARCH model and accommodates conditional heteroskedasticity. In the ARCH/GARCH model, the volatility process is serially dependent and depends on past information. Such dependence is excluded from the elliptical distribution;  however, it is perfectly compatible with  Model~\eqref{observe}.
	\item The dependence between $\{\omega_i\}$ and $\{\fZ_i\}$ can feature the leverage effect in financial econometrics,  {which accounts for the negative correlation between asset return and change in volatility.} Various research has been conducted to study the leverage effect; see, e.g., \cite{schwert1989does}, \cite{campbell1992no}, \cite{ait2013leverage}, \cite{WM14} and \cite{KX17}.
	\item Furthermore, it can capture  (conditional) asymmetry by allowing the entries of $\fZ_i$'s to be asymmetrically distributed.
	The asymmetry is another  stylized fact of financial data.
	For instance, the empirical study in \cite{singleton1986skewness} shows high skewness in individual stock returns. Skewness is also reported in  exchange rate returns in \cite{peiro1999skewness}. \cite{Christoffersen12} documents that asymmetry exists in  standardized returns; see Chapter 6 therein.
\end{itemize}

Because of the heteroskedasticity induced by $\{\omega_i\}$, in this paper we  focus on testing
\[
	H_{0}:~\boldsymbol{\Sigma} \propto \boldsymbol{\Sigma}_0\quad vs. \quad H_{a}:~\boldsymbol{\Sigma} \not\propto\boldsymbol{\Sigma}_0,
\]
in the high-dimensional setting where both $p$ and~$n$ grow to infinity with the ratio $p/n\rightarrow y\in(0,\infty)$.

Note that in many applications, knowing the covariance matrix up to a constant is good enough. For example, the minimum variance portfolio is given by $\boldsymbol{\Sigma}^{-1} \mathbf{1}/(\mathbf{1}^{T}\fSigma^{-1}\mathbf{1})$, where $\mathbf{1} = (1,\ldots,1)^{T}$. The portfolio is therefore invariant to scaling in the covariance matrix.

\subsection{Summary of main results} To deal with heteroskedasticity, we propose to self-normalize the observations. To be specific, we focus on the self-normalized observations $\fY_i/\left|\fY_i\right|$, where $|\cdot|$ stands for the Euclidean norm. Observe that
$$
\frac{\fY_i}{|\fY_i|}=\frac{\fSigma^{1/2}\fZ_i}{|\fSigma^{1/2}\fZ_i|},\quad i=1, \cdots, n.
$$
Hence $\omega_i$'s no longer play a role, and this is exactly the reason why we make no assumption on~$\omega_i$'s. There is, however, no such thing as a free lunch. Self-normalization introduces a new challenge in that the entries of $\fSigma^{1/2}\fZ_i/|\fSigma^{1/2}\fZ_i|$ are dependent
in an unusual fashion.
To see this, consider the simplest case where $\fSigma=\fI$ and $\fZ_i$'s are \hbox{i.i.d.} standard multivariate normal random vectors. In this case, the entries of~$\fZ_i$'s are \hbox{i.i.d.} standard  normal  random variables. However, the self-normalized random vector $\fZ_i/|\fZ_i|$ is uniformly distributed over the $p$-dimensional unit sphere, and its $p$ entries are dependent on each other in an unconventional way.

To conduct tests, we need some kind of CLTs. Our strategy is to establish a CLT for the {\it linear spectral statistic} (LSS)
of the sample covariance matrix based on the self-normalized observations, namely,
\begin{equation}\label{eq:fBn}
\widetilde{\fS}_n=\frac{\tr(\fSigma)}{n}\sum_{i=1}^{n}\frac{\fY_i\fY_i^{T}}{|\fY_i|^2}=\frac{\tr(\fSigma)}{n}\sum_{i=1}^{n}\frac{\fSigma^{1/2}\fZ_i\fZ_i^{\top}\fSigma^{1/2}}{\big|\fSigma^{1/2}\fZ_i\big|^2}.
\end{equation}
When $|\fY_i|$ or $|\fSigma^{1/2}\fZ_i|=0$, we adopt the convention that $0/0=0$.

As we shall see below,  our CLT is different from the ones for the usual sample covariance matrix.  {One important advantage of our result is that, when testing $H_0: \fSigma\propto\fI$, applying our CLT requires neither $\mathbb{E}(Z_{11}^4)=3$ as in \cite{BaiS04}, nor the estimation of $\mathbb{E}(Z_{11}^4)$, which is inevitable in \cite{Najim2013}. {Based on the new CLT, we propose two sphericity tests by modifying the likelihood ratio test and John's test. More tests based on general moments of the ESD of $\widetilde{\fS}_n$ are also constructed.} Numerical studies show that, for the sphericity hypothesis, our proposed tests work well even when $\mathbb{E}(Z_{11}^4)$ does not exist. Because heavy-tailedness and heteroskedasticity are commonly encountered in practice, such relaxations are appealing in many real applications. {More generally, we can also test $H_0: \fSigma\propto\fSigma_0$, where $\fSigma_0$ is a general non-negative definite matrix and can even be not invertible. We illustrate such a test in Section \ref{simu_gen_sig} for a case when $\fSigma_0$ contains a substantial proportion (1/4 to be precise) of zero eigenvalues.} }

{
\begin{rmk}
Independently, \cite{li2017structure}  study high-dimensional covariance matrix test under a mixture model. {Their test relies on comparing two John's test statistics: one is based on the original data and the other is based on randomly permutated data.}
There are a couple of major differences between our paper and theirs. Firstly, they only consider sphericity test and so can only test $H_0: \fSigma\propto\fSigma_0$ when $\fSigma_0$ is invertible. Secondly, in \cite{li2017structure}, the mixture coefficients ($\omega_i$'s in \eqref{observe}) are assumed to be \hbox{i.i.d.} and drawn from a distribution with a bounded support. Thirdly, \cite{li2017structure} require  independence between the mixture coefficients and the innovation process $(\fZ_i)$.
In our paper, we do not put any assumptions on the mixture coefficients. As we discussed in Section \ref{Ourmodelandgoal}, such relaxations allow us to accommodate several important stylized features of real data,  consequently, make our  tests more suitable in many real applications.
It can be shown that the test in \cite{li2017structure} is not necessarily consistent under our setup even if $\fSigma$ in our model \eqref{observe} is identity. Furthermore, as we can see from the simulation studies,  the test in \cite{li2017structure} is less powerful than the existing tests {in the \hbox{i.i.d.} Gaussian setting and}, in general, substantially less powerful than our tests.
\end{rmk}}

Empirically, we apply the proposed tests to study the correlations among idiosyncratic returns. There have been numerous factor models developed. A very interesting question is whether the idiosyncratic returns under a factor model are uncorrelated. If the answer is yes, then one can conclude that there is no missing factor. Answering such a question is challenging because typically there are tens or  hundreds of stocks involved, resulting in thousands of correlations to be tested.
As an innovative application, we demonstrate that  our tests can be utilized to test uncorrelatedness among idiosyncratic returns. We illustrate the testing procedure by using the CAPM (\cite{Sharpe64} and the Fama-French Three-Factor model (\cite{FF92}). The analysis can be translated directly to more comprehensive factor models.

\bigskip

The rest of the paper is organized as follows. In Section~\ref{MRTS}, we state the CLT for the LSS of~$\widetilde{\fS}_n$, based on which, for testing $H_0: \fSigma\propto\fI$, we derive the asymptotic null distributions of the modified likelihood ratio test statistic and John's test statistic, {as well as other test statistics based on general moments of the ESD of $\widetilde{\fS}_n$.} {Section~\ref{testpro} examines the finite-sample performance of the proposed sphericity tests, and illustrates how to test $H_0: \fSigma\propto\fSigma_0$ for a non-invertible non-negative definite matrix $\fSigma_0$.}  {Section~\ref{sect:empirical} is dedicated to a real data analysis, in which we show how our tests can be used to test uncorrelatedness among idiosyncratic returns.}
Section~\ref{conclusion} concludes. All  proofs are collected in Appendices \ref{proof_Thm_CLTLSS}--\ref{two_lemma}.

Finally, we collect some notation that will be used throughout the paper. For any symmetric matrix $\fA\in\mathbb{R}^{p\times p}$, $\|\fA\|$ stands for the spectral norm and $F^{\fA}$ denotes the ESD, that is,
$$
\|\fA\|=\max_i|\lambda^{\fA}_i|, \text{ and }F^{\fA}(x)=\frac{1}{p}\sum_{i=1}^p \mathbbm{1}_{\{\lambda_i^{\fA}\leq x\}},~~ \mbox{for all } x\in\mathbb{R},
$$
where $\lambda^{\fA}_i$, $i=1, \cdots, p$, are the eigenvalues of $\fA$ and $\mathbbm{1}_{\{\cdot\}}$ denotes the indicator function. We also denote by $\lambda_{\min}^\fA$ the smallest eigenvalue of $\fA$.
For any function~$f$, the associated LSS of $\fA$ is given by
\begin{align*}
\int_{-\infty}^{+\infty}f(x){\rm d}F^{\fA}(x)=\frac{1}{p}\sum_{i=1}^pf(\lambda_i^{\fA}).
\end{align*}
The Stieltjes transform of a distribution $G$ is defined as
\begin{align*}
m_G(z)=\int_{-\infty}^{\infty}\frac{1}{\lambda-z}{\rm d}G(\lambda), \quad\mbox{for all } z\not\in\hbox{supp}(G),
\end{align*}
where supp($G$) denotes the support of $G$. For any $y\in(0, \infty)$ and distribution $G$, $F_{y, G}$ denotes the distribution whose Stieltjes transform $m_{F_{y, G}}(z)$ is the unique solution to
\[
	m(z)=\int_0^{\infty}\frac{{\rm d}G(t)}{t(1-y-yzm(z))-z}, ~z\in\mathbb{C}^+:=\{z\in\mathbb{C}, \Im(z)>0\},
\]	
in the set $\{m(z)\in\mathbb{C}: -(1-y)/z+ym(z)\in\mathbb{C}^+\}$, where $\Im(z)$ denotes the imaginary part of~$z$. Finally, for any $x\in\mathbb{R}$, $\delta_{x}$ is the Dirac measure at~$x$, and we denote $F_{y, \delta_1}$ as $F_{y}$.

\section{Main Results}\label{MRTS}

\subsection{CLT for the LSS of $\widetilde{\fS}_n$}\label{ssec:CLT_LSS}
As discussed above, we focus on the sample covariance matrix based on the self-normalized observations, namely,
$\widetilde{\fS}_n$ defined in~\eqref{eq:fBn}.

We make the following assumptions:

{\begin{myassump}{A}\label{ass:A}
The random variables $Z_{ij}$'s are \mbox{i.i.d.} with $\mathbb{E}(Z_{11})=0$, $\mathbb{E}(Z_{11}^2)=1$  and $\mathbb{E}(Z_{11}^{4+\zeta})<\infty$ for some $\zeta>0$.
\end{myassump}}

{\begin{myassump}{A$^\prime$}\label{ass:A'}
The random variables $Z_{ij}$'s are \mbox{i.i.d.} with $\mathbb{E}(Z_{11})=0$, $\mathbb{E}(Z_{11}^2)=1$  and $\mathbb{E}(Z_{11}^{4})<\infty$.
\end{myassump}}

{
\begin{myassump}{B}\label{ass:B}
The probability density function, $f_Z(\cdot)$, of $Z_{11}$ satisfies $0\leq f_Z(\cdot)\leq C_1$ for some $C_1>0$.
\end{myassump}}

\begin{myassump}{C}\label{ass:C}
There exists a distribution $H$ such that
$H_p\stackrel{D}\longrightarrow H$ as $p\rightarrow\infty$, where $H_p=F^{\fSigma}$. Furthermore, $\tr\big(\fSigma\big)\asymp p$,  $\|\fSigma\|\leq C_2$ and $\widetilde{\lambda}^{\fSigma}_{\min}\geq p^{-C_3}$ for some constants $C_2>0$ and $C_3>0$, where $\widetilde{\lambda}^{\fSigma}_{\min}$ denotes the smallest non-zero eigenvalue of $\fSigma$; and
\end{myassump}
\begin{myassump}{D}\label{ass:D} $y_n:=p/n\rightarrow y\in\left(0, \infty\right)$ as $n\rightarrow\infty$.
\end{myassump}

Theorem 2 in \cite{ZhengL11} states that under some regularity conditions, $\widetilde{\fS}_n$ shares the same LSD as the sample covariance matrix~$\fS_n:=n^{-1}\sum_{i=1}^n\fSigma^{1/2}\fZ_i\fZ_i^\top\fSigma^{1/2}$. To conduct tests, we need the associated CLT. The CLTs for the LSS of~$\fS_n$ have been established in \cite{BaiS04} and \cite{Najim2013}, {under the Gaussian and {non-Gaussian kurtosis conditions}, respectively}. Given that $\widetilde{\fS}_n$ and~$\fS_n$ have the same LSD, one naturally asks  whether their LSSs also have the same CLT. The following theorem gives a negative answer. Hence, an important message is:
\begin{center}
	\emph{Self-normalization does not change the LSD, but it does affect the CLT.}
\end{center}

To be more specific, for any function $f$, define the following centered and scaled LSS:
\begin{align}\label{generalG}
G_{\widetilde{\fS}_n}(f):=p\int_{0}^{+\infty}f(x)\ {\rm d} \Big(F^{\widetilde{\fS}_n}(x)-F_{y_n, H_p}(x)\Big).
\end{align}
{\begin{thm}\label{CLTLSS}
Let $a_l(\fSigma, y)=\ul{\lim}_{n\rightarrow\infty}\lambda_{\min}^{\fSigma}\mathbbm{1}_{\{0<y<1\}}(1-\sqrt{y})^2,$ $a_r(\fSigma, y)=\overline{\lim}_{n\rightarrow\infty}\|\fSigma\|(1+\sqrt{y})^2$, $\mathcal{H}$ be the set of functions that are analytic on a domain containing $[a_l(\fSigma, y), a_r(\fSigma, y)]$, and $f_1, \ldots, f_k\in\mathcal{H}$.
	\begin{enumerate}[(i)]
\item\label{item:Meq3}
Under Assumptions \ref{ass:A}--\ref{ass:D}, the sequence of random vectors $\big\{(G_{\widetilde{\fS}_n}(f_1)$, $\ldots$, $G_{\widetilde{\fS}_n}(f_k))\big\}$ is tight. Furthermore, if $\mathbb{E}(Z_{11}^4)=3$, then the random vector $\big(G_{\widetilde{\fS}_n}(f_1)$, $\ldots$, $G_{\widetilde{\fS}_n}(f_k)\big)$ converges weakly to a Gaussian vector $\big(G(f_1),$ $\ldots, $ $G(f_k)\big)$ with mean
\begin{equation}
\hspace{-1em}\begin{aligned}\label{eq:EGf_first}
	\mathbb{E}\big(G(f_i)\big)
	=&-\frac{1}{2\pi\mathrm{i}}\oint_{\mathcal{C}} f_i(z)\int_0^\infty\frac{y\ul{m}^3(z)t^2{\rm d}H(t)}{(1+t\ul{m}(z))^3}\bigg(\!1\!-\!\!\int_0^\infty\frac{y\ul{m}^2(z)t^2{\rm d}H(t)}{(1+t\ul{m}(z))^2}\bigg)^{\!-2}\!\!{\rm d}z\\
	&+\frac{1}{2\pi\mathrm{i}}\oint_{\mathcal{C}} f_i(z)R_0(z)\ul{m}(z)\bigg(1-\int_0^\infty\frac{y\ul{m}^2(z)t^2{\rm d}H(t)}{(1+t\ul{m}(z))^2}\bigg)^{-1}{\rm d}z,
	~~i=1, \ldots, k,
		\end{aligned}
	\end{equation}
	where
	\begin{equation*}
\begin{aligned}
R_0(z)=&\frac{2\int_0^\infty t^2{\rm d}H(t)}{(\int_0^\infty t{\rm d}H(t))^2}\int_0^\infty\frac{t\ul{m}(z){\rm d}H(t)}{(1+t\ul{m}(z))^2}
-\frac{2}{\int_0^\infty t{\rm d}H(t)}\int_0^\infty\frac{t^2\ul{m}(z){\rm d}H(t)}{(1+t\ul{m}(z))^2}\\
&+\frac{2y}{\int_0^\infty t{\rm d}H(t)}\int_0^\infty\frac{t\ul{m}(z){\rm d}H(t)}{1+t\ul{m}(z)}\int_0^\infty\frac{t^2\ul{m}(z){\rm d}H(t)}{(1+t\ul{m}(z))^2}\\
&+\frac{2y}{\int_0^\infty t{\rm d}H(t)}\int_0^\infty\frac{t^2\ul{m}(z){\rm d}H(t)}{1+t\ul{m}(z)}\int_0^\infty\frac{t\ul{m}(z){\rm d}H(t)}{(1+t\ul{m}(z))^2}\\
&-\frac{2y\int_0^\infty t^2{\rm d}H(t)}{(\int_0^\infty t{\rm d}H(t))^2}\int_0^\infty\frac{t\ul{m}(z){\rm d}H(t)}{1+t\ul{m}(z)}\int_0^\infty\frac{t\ul{m}(z){\rm d}H(t)}{(1+t\ul{m}(z))^2},
\end{aligned}
\end{equation*}
and covariance
\begin{equation}
	\begin{aligned}\label{eq:thmcov_first}
	&\cov((G(f_i), G(f_j))\\
	=&-\frac{1}{2\pi^2}\oint_{\mathcal{C}_2}\!\!\oint_{\mathcal{C}_1}\!\!\frac{f_i(z_1)f_j(z_2)\ul{m}^{\prime}(z_1)\ul{m}^{\prime}(z_2)}{\big(\ul{m}(z_2)-\ul{m}(z_1)\big)^2}{\rm d}z_1{\rm d}z_2\\
	&+\frac{y}{2\pi^2\int_0^\infty\!\!t{\rm d}H(t)}\oint_{\mathcal{C}_2}\!\!\oint_{\mathcal{C}_1}\!\!\int_0^\infty\!\!\frac{tf_i(z_1)\ul{m}^{\prime}(z_1){\rm d}H(t)}{(1+t\ul{m}(z_1))^2}\!\!\int_0^\infty\!\!\frac{t^2f_j(z_2)\ul{m}^{\prime}(z_2){\rm d}H(t)}{(1+t\ul{m}(z_2))^2}{\rm d}z_1{\rm d}z_2\\
	&+\frac{y}{2\pi^2\int_0^\infty t{\rm d}H(t)}\oint_{\mathcal{C}_2}\!\!\oint_{\mathcal{C}_1}\!\!\int_0^\infty\!\!\frac{t^2f_i(z_1)\ul{m}^{\prime}(z_1){\rm d}H(t)}{(1+t\ul{m}(z_1))^2}\!\!\int_0^\infty\!\!\frac{tf_j(z_2)\ul{m}^{\prime}(z_2){\rm d}H(t)}{(1+t\ul{m}(z_2))^2}{\rm d}z_1{\rm d}z_2\\
	&-\frac{y\int_0^\infty\!\!t^2{\rm d}H(t)}{2\pi^2\big(\int_0^\infty\!\!t{\rm d}H(t)\big)^2}\oint_{\mathcal{C}_2}\!\!\oint_{\mathcal{C}_1}\!\!\int_0^\infty\!\!\frac{tf_i(z_1)\ul{m}^{\prime}(z_1){\rm d}H(t)}{(1+t\ul{m}(z_1))^2}\!\!\int_0^\infty\frac{tf_j(z_2)\ul{m}^{\prime}(z_2){\rm d}H(t)}{(1+t\ul{m}(z_2))^2}{\rm d}z_1{\rm d}z_2,
	\end{aligned}
	\end{equation}
	$i, j=1, \ldots, k$, where
	$\ul{m}(z)$ is the Stieltjes transform of $\underline{F}_{y, H}:=(1-y)\mathbbm{1}_{[0,\infty)}+yF_{y, H}$, and
	$\mathcal{C}_1$ and $\mathcal{C}_2$ are two non-overlapping contours contained in the domain and enclosing the interval $[a_l(\fSigma, y), a_r(\fSigma, y)]$;
	\item If $\fSigma=\fI$, then under Assumptions \ref{ass:A'} and \ref{ass:D} and without assuming $\mathbb{E}(Z_{11}^4)=3$, the weak convergence in (i) still holds, and the mean and covariance admit the following simpler expressions:
	\begin{equation}
	\begin{aligned}\label{eq:EGf}
	\hspace{-1em}\mathbb{E}\big(G(f_i)\big)
	=&\!-\!\frac{1}{2\pi\mathrm{i}}\oint_{\mathcal{C}} f_i(z)\Bigg(\frac{y\ul{m}^3(z)}{\big(1+\ul{m}(z)\big)^3}\Bigg)\!\Bigg(1-\frac{y\ul{m}^2(z)}{\big(1+\ul{m}(z)\big)^2}\Bigg)^{-2}{\rm d}z\\
	&+\frac{1}{\pi\mathrm{i}}\oint_{\mathcal{C}} f_i(z)\Bigg(\frac{y\ul{m}^3(z)}{\big(1+\ul{m}(z)\big)^3}\Bigg)\!\Bigg(1-\frac{y\ul{m}^2(z)}{\big(1+\ul{m}(z)\big)^2}\Bigg)^{-1}{\rm d}z;
	\end{aligned}
	\end{equation}
	\begin{equation}
	\begin{aligned}\label{eq:thmcov}
	\cov((G(f_i), G(f_j))=&-\frac{1}{2\pi^2}\oint_{\mathcal{C}_2}\oint_{\mathcal{C}_1}\frac{f_i(z_1)f_j(z_2)\ul{m}'(z_1)\ul{m}'(z_2)}{\big(\ul{m}(z_2)-\ul{m}(z_1)\big)^2}{\rm d}z_1{\rm d}z_2\\
	&+\frac{y}{2\pi^2}\oint_{\mathcal{C}_2}\oint_{\mathcal{C}_1}\frac{f_i(z_1)f_j(z_2)\ul{m}^{\prime}(z_1)\ul{m}^{\prime}(z_2)}{\big(1+\ul{m}(z_1)\big)^2\big(1+\ul{m}(z_2)\big)^2}{\rm d}z_1{\rm d}z_2,
	\end{aligned}
	\end{equation}
	where $i, j=1, \ldots, k.$
\end{enumerate}	

\end{thm}}



{
\begin{rmk}
The first terms in \eqref{eq:EGf_first} and \eqref{eq:thmcov_first} appear in equations~(1.6) and (1.7) of \cite{BaiS04}.
The first terms in \eqref{eq:EGf} and \eqref{eq:thmcov} equal the terms in equations~(1.6) and (1.7) of \cite{BaiS04} with $\fSigma=\fI$.
The other terms in~\eqref{eq:EGf_first}--\eqref{eq:thmcov} are new and are due to the self-normalization in $\widetilde{\fS}_n$. It is worth emphasizing that, when $\fSigma=\fI$,
our CLT neither requires $\mathbb{E}(Z_{11}^4)=3$ as in \cite{BaiS04}, nor involves $\mathbb{E}(Z_{11}^4)$ as in \cite{Najim2013}.
\end{rmk}}

{
\begin{rmk}
Our CLT allows $\fSigma$ to be not invertible and so can be used to test $H_0: \fSigma\propto\fSigma_0$ even when $\fSigma_0$ is not invertible. This is an important contribution in the covariance matrix testing literature because existing methods rely on transforming testing $H_0: \fSigma\propto\fSigma_0$ to $H_0: \fSigma\propto\fI$ by multiplying observations with $\fSigma_0^{-1/2}$.
\end{rmk}
}

\subsection{Tests of sphericity in the presence of heteroskedasticity} \label{prop_test_fi}
For sphericity test
\begin{align}\label{test:identityI}
H_{0}:~\boldsymbol{\Sigma} \propto \fI\quad vs. \quad H_{a}:~\boldsymbol{\Sigma} \not\propto\fI,
\end{align}
based on Theorem~\ref{CLTLSS}, we propose two tests by modifying the likelihood ratio test and John's test. {More tests based on general moments of the ESD of $\widetilde{\fS}_n$  are also established.}

\subsubsection{Likelihood ratio test based on self-normalized observations (LR-SN)}

The classical likelihood ratio test statistic is
$$
L_n=\log|\fS_n|-p\log\big(\tr\big(\fS_n\big)\big)+p\log p;
$$
see, e.g.,  Section 8.3.1 in \cite{Muirhead82}.
For the heteroskedastic case, we modify the likelihood ratio test statistic by replacing $\fS_n$ with $\widetilde{\fS}_n$. Note that $\tr\big(\widetilde{\fS}_n\big)=p$ on the event $\{|\fZ_i|>0~\text{for}~ i=1, \ldots, n\}$, which, by Lemma 2 in \cite{BaiY93}, occurs almost surely for all large $n$. Therefore, we are led to the following modified likelihood ratio test statistic:
\begin{align*}
\widetilde{L}_n=\log\big|\widetilde{\fS}_n\big|=\sum_{i=1}^p\log\big(\lambda^{\widetilde{\fS}_n}_i\big).
\end{align*}
It is the LSS of $\widetilde{\fS}_n$ when $f(x)=\log(x)$.  In this case, when $y_n\in (0,1),$ we have
\begin{align*}
G_{\widetilde{\fS}_n}(\log)=&p\int_{-\infty}^{+\infty}\log(x){\rm d}\left(F^{\widetilde{\fS}_n}(x)-F_{y_n}(x)\right)\\
=&\sum_{i=1}^p\log\big(\lambda^{\widetilde{\fS}_n}_i\big)-p\bigg(\frac{y_n-1}{y_n}\log(1-y_n)-1\bigg)\\
=&\widetilde{L}_n-p\bigg(\frac{y_n-1}{y_n}\log(1-y_n)-1\bigg).
\end{align*}

Applying Theorem \ref{CLTLSS}, we obtain the following proposition.
\begin{prop}\label{prop:CLT_specific_log}
When $y_n\rightarrow y\in(0, 1)$, under Assumption \ref{ass:A'}, we have
	\begin{align}
	\frac{\widetilde{L}_n-p\bigg(\frac{y_n-1}{y_n}\log(1-y_n)-1\bigg)-\big(\log(1-y_n)\big)/2-y_n}{\sqrt{-2\log(1-y_n)-2y_n}}\stackrel{D}\longrightarrow N(0, 1). \label{eq:CLT_log}
	\end{align}
\end{prop}

The convergence in \eqref{eq:CLT_log} gives the asymptotic null distribution of the modified likelihood ratio test statistic. Because it is derived for the sample covariance matrix based on self-normalized observations, the test based on \eqref{eq:CLT_log} will be referred to as the likelihood ratio test based on the self-normalized observations (LR-SN).

\subsubsection{John's test based on self-normalized observations (JHN-SN)}

John's test statistic (\cite{John1971}) is given by
$$
T_n=\frac{n}{p}\tr\bigg(\frac{\fS_n}{1/p\tr\big(\fS_n\big)}-\fI\bigg)^2-p.
$$
Replacing $\fS_n$ with $\widetilde{\fS}_n$ and noting again that $\tr\big(\widetilde{\fS}_n\big)=p$ almost surely for all large~$n$ lead to the following modified John's test statistic:
\begin{align*}
\widetilde{T}_n=\frac{n}{p}\tr\Big(\widetilde{\fS}_n-\fI\Big)^2
-p=\frac{1}{y_n}\sum_{i=1}^p\big(\lambda^{\widetilde{\fS}_n}_i\big)^2-n-p.
\end{align*}
It is related to the LSS of $\widetilde{\fS}_n$ when $f(x)=x^2$. In this case, we have
\begin{align*}
G_{\widetilde{\fS}_n}(x^2)=p\int_{-\infty}^{+\infty}x^2\ {\rm d}\left(F^{\widetilde{\fS}_n}(x)\!-\!F_{y_n}(x)\right)=\sum_{i=1}^p\big(\lambda^{\widetilde{\fS}_n}_i\big)^2-p(1+y_n)
=y_n\widetilde{T}_n.
\end{align*}

Based on Theorem \ref{CLTLSS}, we can prove the following proposition.
\begin{prop}\label{prop:CLT_specific_x^2}
Under Assumptions \ref{ass:A'} and \ref{ass:D},  we have
	\begin{align}
	\frac{\widetilde{T}_n+1}{2}\stackrel{D}\longrightarrow N(0, 1).\label{eq:CLT_sq}
	\end{align}
\end{prop}

Below we will refer to the test based on \eqref{eq:CLT_sq} as John's test based on the self-normalized observations (JHN-SN).

\subsubsection{More general tests based on self-normalized observations }

More tests can be constructed  by choosing $f$ in Theorem \ref{CLTLSS} to be different functions. When $f(x)=x^k$ for $k\geq2$, the corresponding LSS is the
$k$th moment of the ESD of $\widetilde{\fS}_n$, for which we have
	\begin{align*}
	G_{\widetilde{\fS}_n}(x^k)=&p\int_{-\infty}^{+\infty}x^k\ {\rm d}\left(F^{\widetilde{\fS}_n}(x)\!-\!F_{y_n}(x)\right)\\
	=&\sum_{i=1}^p\big(\lambda^{\widetilde{\fS}_n}_i\big)^k-p(1+y_n)^{k-1}H_F\Big(\frac{1-k}{2}, 1-\frac{k}{2}, 2, \frac{4y_n}{(1+y_n)^2}\Big),
	\end{align*}
	where $H_F(a, b, c, d)$ denotes the hypergeometric function $_2F_1(a, b, c, d)$. By Theorem \ref{CLTLSS} again, we have the following proposition.
	\begin{prop}\label{prop:CLT_x^k}
	Under Assumptions \ref{ass:A'} and \ref{ass:D}, for any $k\geq2$,  we have
		\begin{align*}
		\frac{G_{\widetilde{\fS}_n}(x^k)-\mu_{n, x^k}}{\sigma_{n, x^k}}\stackrel{D}\longrightarrow N(0, 1),\quad \text{where}
		\end{align*}
		\begin{align*}
		\mu_{n, x^k}=&\frac{1}{4}\Big((1+\sqrt{y_n})^{2k}+(1-\sqrt{y_n})^{2k}\Big)-\frac{1}{2}\sum_{i=0}^k{k \choose i}^2y_n^i\\
		&-\frac{2k(k-1)(1+y_n)^{k-2}}{(k+1)(k+2)}\bigg(\!(y_n-1)^2H_F\Big(\!\frac{3-k}{2}, 1-\frac{k}{2}, 1, \frac{4y_n}{(1+y_n)^2}\!\Big)\\
		&~~~~~~~~+(-1+4ky_n-y_n^2)H_F\Big(\frac{3-k}{2}, 1-\frac{k}{2}, 2, \frac{4y_n}{(1+y_n)^2}\Big)\!\bigg),
\end{align*}
and
\begin{align*}
		\sigma^2_{n, x^k}=&2y_n^{2k}\sum_{i=0}^{k-1}\sum_{j=0}^k{k\choose i}{k\choose j}\Big(\frac{1-y_n}{y_n}\Big)^{i+j}\sum_{\ell=1}^{k-i}\ell{2k-1-(i+\ell)\choose k-1}{2k-1-j+\ell\choose k-1}\\
		&-2y_n\bigg((1-y_n)^kk\sum_{i=0}^{k+1}{k+1\choose i}\Big(\frac{1-y_n}{y_n}\Big)^{1-i}\frac{(k+i-1)!}{(i-1)!(k+1)!}\bigg)^2.
	\end{align*}
	\end{prop}

\begin{rmk}
	Proposition \ref{prop:CLT_x^k} enables us to consistently detect any alternative hypothesis under which the covariance matrix admits an LSD not equal to {$\delta_{1}$}. The reason is that, under such a situation, the LSD of $\widetilde{\fS}_n$, say $\widetilde{H}$, will not be the standard Mar\v{c}enko-Pastur law $F_y$. Therefore,  there exists a $k\geq 2$  such that $\int_{-\infty}^{\infty} x^k\, {\rm d}\widetilde{H}(x)\neq\int_{-\infty}^{\infty} x^k\, {\rm d}F_y(x)$. Consequently,  $G_{\widetilde{\fS}_n}(x^k)$ in \eqref{generalG} will blow up, and the testing power will approach $1$.
\end{rmk}

\section{Simulation Studies}\label{testpro}

We now demonstrate the finite-sample performance of our proposed tests in Section \ref{prop_test_fi}, as well as the test for
$H_{0}: \boldsymbol{\Sigma} \propto \boldsymbol{\Sigma}_0$, where $\fSigma_0$ is a general non-negative definite matrix.

\subsection{\hbox{I.i.d.} Gaussian case} To have a full picture of the performance of our tests, we start with the simplest situation where observations are \hbox{i.i.d.} multivariate normal random vectors.
{We will compare our proposed tests,  LR-SN and JHN-SN, with the tests mentioned in Section \ref{test_hd}, namely,  LW$_2$, S, CZZ$_2$ and WY-LR, and also the test in \cite{li2017structure} (LY test).}  In the multivariate normal case, the WY-JHN test reduces to the LW$_2$ test.

We start with the size evaluation by sampling  observations  from $N(\mathbf{0}, \fI)$. Table \ref{compare_size} reports the empirical sizes of these tests for testing $H_0: \boldsymbol{\Sigma}\propto\fI$  at $5\%$ significance level.

\begin{center}
+++  Insert Table \ref{compare_size} Here +++
\end{center}

From Table \ref{compare_size}, we find that the empirical sizes of all tests are around the nominal level of $5\%$.

Next, to compare the power, we
generate \hbox{i.i.d.} observations from $N(\mathbf{0}, \boldsymbol{\Sigma})$ under the alternative with $\boldsymbol{\Sigma}=\big(0.1^{|i-j|}\big)_{p\times p}$, and test $H_{0}: \boldsymbol{\Sigma}\propto\fI$
 at $5\%$ significance level. Table \ref{compare_power} reports the empirical powers.

\begin{center}
+++  Insert Table \ref{compare_power} Here +++
\end{center}

From Table \ref{compare_power}, we find that our proposed LR-SN and JHN-SN  tests and the tests mentioned in Section \ref{test_hd} have quite high powers especially as the dimension gets higher, and the powers are roughly comparable. Same as in the classical setting, John's test (JHN-SN) is more powerful than the likelihood ratio test (LR-SN). LY test proposed in \cite{li2017structure} is less powerful.

{To sum up, while developed under a much more general setup, our tests perform just as well as the existing ones in the ideal \hbox{i.i.d.} Gaussian setting.}

\subsection{The elliptical case}  We now investigate the performance of our proposed tests under the elliptical distribution.  As in Section \ref{performance_of_existing_tests}, we take the observations to be $\fY_i=\omega_i\fZ_i$ with
\begin{enumerate}[(i)]
	\item $\omega_i$'s being absolute values of \hbox{i.i.d.} standard normal random variables,
	\item $\fZ_i$'s \hbox{i.i.d.} $p$-dimensional random vectors from $N(\mathbf{0}, \boldsymbol{\Sigma})$, and
	\item $\omega_i$'s and $\fZ_i$'s independent of each other.
\end{enumerate}

\emph{\textbf{Checking the size.}}

Table \ref{ellipticalsize} completes Table \ref{counterexample_sphere} by including the empirical sizes of our proposed LR-SN and JHN-SN tests, and also LY test in \cite{li2017structure}.

\begin{center}
+++  Insert Table \ref{ellipticalsize} Here +++
\end{center}

Table \ref{ellipticalsize} reveals sharp difference between the existing tests and our proposed ones: the empirical sizes of the existing tests are severely distorted, in contrast,  the empirical sizes of our LR-SN and JHN-SN tests are around the nominal level of $5\%$ as desired. LY test also yields the right level of size.

\bigskip

\emph{\textbf{Checking the power.}}

Table \ref{ellipticalsize} shows that LW$_2$, S, CZZ$_2$, WY-LR and WY-JHN tests are inconsistent under the elliptical distribution, therefore we exclude them when checking the power.

We generate observations  under the elliptical distribution with $\boldsymbol{\Sigma}=\big(0.1^{|i-j|}\big)$. Table~\ref{ellipticalpower} reports the empirical powers of our proposed tests and LY test for testing $H_{0}: \boldsymbol{\Sigma}\propto\fI$ at $5\%$ significance level.

\begin{center}
+++  Insert Table \ref{ellipticalpower} Here +++
\end{center}

From Table \ref{ellipticalpower}, we find that
\begin{enumerate}[(i)]
  \item Our tests, LR-SN and JHN-SN, as well as LY test, enjoy a blessing of dimensionality: for a fixed ratio $p/n$, the higher the dimension $p$, the higher the power;
  \item LY test is substantially less powerful than our tests.
\end{enumerate}

\subsection{Beyond elliptical, a GARCH-type case}\label{ssec:sim_hetero_case}
Recall that in our general model \eqref{observe}, the observations $\fY_i$ admit the decomposition $\omega_i\fSigma^{1/2}\fZ_i$, and $\omega_i$'s can depend on each other and on $\{\fZ_i: i=1, \ldots, n\}$ in an arbitrary way. To examine the performance of our tests in such a general setup, we simulate data using the following two-step procedure:
\begin{enumerate}[1.]
{\item Sample ${Z}_{ij}$'s from standardized $t$-distribution with four degrees of freedom. The distribution is heavy-tailed and even does not have finite fourth moment.}
	\item For each $\omega_i$, inspired by the ARCH/GARCH model, we take $\omega^2_i=0.01+0.85\omega^2_{i-1}+0.1|\fY_{i-1}|^2/\tr\big(\boldsymbol{\Sigma}\big)$.
\end{enumerate}

\emph{\textbf{Checking the size.}}

We take $\fSigma=\fI$ in the data generating process and test $H_{0}: \boldsymbol{\Sigma}\propto\fI$.
Table~\ref{size} reports the empirical sizes of our proposed tests and LY test at $5\%$ significance level.

\begin{center}
+++  Insert Table \ref{size} Here +++
\end{center}

From Table \ref{size}, we find that, for all different values of $p$ and $p/n$, the empirical sizes of our proposed tests are around the nominal level of $5\%$. Again, this is in sharp contrast with the results in Table \ref{counterexample_sphere}, where the existing tests yield sizes far higher than $5\%$.

The second observation is that although Theorem \ref{CLTLSS} requires the finiteness of~$\mathbb{E}\big(Z_{11}^{4}\big)$, the simulation above shows that
our proposed tests work well even when $\mathbb{E}(Z_{11}^4)$ does not exist.

Another observation is that the sizes of LY test are around $8\%$. Note that with 10,000 replications, the margin of error at $5\%$ significance level is $1\%$, hence we can conclude that the sizes of LY test are statistically significantly higher than the nominal level of $5\%$.

\bigskip
\emph{\textbf{Checking the power.}}

To evaluate the power, we again take  $\boldsymbol{\Sigma}=\big(0.1^{|i-j|}\big)$ and generate data according to the design at the beginning of this subsection. Table \ref{tab:power} reports the empirical powers of our proposed tests and LY test for testing $H_{0}: \boldsymbol{\Sigma}\propto\fI$ at $5\%$ significance level.

\begin{center}
+++  Insert Table \ref{tab:power} Here +++
\end{center}

{Table \ref{tab:power} shows again that our tests enjoy a blessing of dimensionality. Moreover, comparing Table \ref{tab:power} with Table \ref{ellipticalpower}, we find that for each pair of $p$ and $n$, the powers of our tests are similar under the two designs. Such similarities show that our tests can not only accommodate (conditional) heteroskedasticity but also are robust to  heavy-tailedness in~$\fZ_i$'s. Finally, LY test is again significantly less powerful.}

{\subsection{A GARCH-type case with general $\fSigma_0$}\label{simu_gen_sig}
In this subsection, we first illustrate how to construct the test for general $\fSigma_0$ based on Theorem \ref{CLTLSS} and
then investigate the finite-sample performance of the test.
We take $\fSigma_0=\fQ\boldsymbol{\Lambda}_0\fQ^\top$, where $\boldsymbol{\Lambda}_0=\text{diag}(\underbrace{1, \ldots, 1}_{p/2}, \underbrace{2, \ldots, 2}_{p/4}, \underbrace{0, \ldots, 0}_{p/4})$ and $\fQ$ is a random orthogonal matrix. To build the test, we take $f(x)=x^2$.

We need to compute the asymptotic mean in \eqref{eq:EGf_first}. Note that by
equation (1.2) in \citeapp{BaiS04}, we have
\begin{align}\label{eq:zeqlm}
z(\ul{m})=-\frac{1}{\ul{m}}+y\Big(\frac{1}{2(1+\ul{m})}+\frac{1}{2(1+2\ul{m})}\Big),
\end{align}
where $\ul{m}=\ul{m}(z)$. By the Chain rule with the transformation \eqref{eq:zeqlm}, we can find that the asymptotic mean  in \eqref{eq:EGf_first} equals the sum of some contour integrals with respective to $\ul{m}$, and the singularities are $-1$, $-1/2$ and the roots of
\begin{align*}
\bigg(1-\int_0^{\infty}\frac{y\ul{m}^2t^2{\rm d}H(t)}{(1+t\ul{m})^2}\bigg)^{-1}=0, \text{where }{H=\frac{1}{2}\delta_{1}+\frac{1}{4}\delta_{2}+\frac{1}{4}\delta_{0}.}
\end{align*}
Then we can compute the asymptotic mean in \eqref{eq:EGf_first} using the residual theorem. Similarly, note that ${\rm d}\ul{m}=\ul{m}^{\prime}(z){\rm d}z$. By the Chain rule, \eqref{eq:zeqlm} and the residue theorem, we can compute the asymptotic covariance in \eqref{eq:thmcov_first}.

The observations are still generated according to the design in Section \ref{ssec:sim_hetero_case} except that ${Z}_{ij}$'s are sampled from the standard normal distribution.

\bigskip
\emph{\textbf{Checking the size.}}

We take $\fSigma=\fSigma_0$ to generate the data and test $H_0:\fSigma\propto\fSigma_0$. Table \ref{tab:size_general} reports the empirical sizes of the test at $5\%$ significance level.

\begin{center}
+++  Insert Table \ref{tab:size_general} Here +++
\end{center}

We find that from Table \ref{tab:size_general} the empirical sizes of the test are around the nominal level of $5\%$ for all different values of $p$ and $p/n$.

\bigskip
\emph{\textbf{Checking the power.}}

To examine the power, we generate the observations with $\fSigma=\fQ\boldsymbol{\Lambda}\fQ^\top$, where $\boldsymbol{\Lambda}=\text{diag}(\underbrace{1, \ldots, 1}_{p/2}, $ $\underbrace{2.5, \ldots, 2.5}_{p/4}, \underbrace{0, \ldots, 0}_{p/4})$ and $\fQ$ is a random orthogonal matrix. Table \ref{tab:power_general} reports the empirical powers of our test for testing $H_0:\fSigma\propto\fSigma_0$ at $5\%$ significance level.

\begin{center}
+++  Insert Table \ref{tab:power_general} Here +++
\end{center}

Again, we see that our test enjoys the blessing of dimensionality and has high powers.  }

\subsection{Summary of simulation studies}\label{ssec:summary_simulation}
Combining the observations in the four cases, we conclude that
\begin{enumerate}[(i)]
	\item {The existing tests, LW$_2$, S, CZZ$_2$, WY-LR and WY-JHN, work well in the \hbox{i.i.d.} Gaussian setting, however, they fail badly under the elliptical distribution and our general setup;}
	\item The newly proposed LY test in \cite{li2017structure} is applicable to the elliptical distribution, however, it is less powerful than the existing tests {in the \hbox{i.i.d.} Gaussian} setting and substantially less powerful than ours in general situations;
	\item Our LR-SN and JHN-SN tests perform well under all the settings, yielding the right sizes and enjoying high powers;
	\item {Even when $\fSigma_0$ is not invertible, our test still works well.}
\end{enumerate}

{\section{Empirical Studies}\label{sect:empirical}}

Let us first explain the motivation of the empirical study, which is about idiosyncratic returns. In general,
the total risk of a stock return can be decomposed into two components: systematic risk and idiosyncratic risk. Empirical studies in \cite{campbell2001have} and \cite{goyal2003idiosyncratic} show that idiosyncratic risk is the major component of the total risk.
It is not uncommon to assume that idiosyncratic returns are cross-sectionally uncorrelated, giving rise to the so-called strict factor model;
see, e.g., \cite{roll1980empirical}, \cite{brown1989number} and \cite{fan2008high}.
Our goal in this section is to test the cross-sectional uncorrelatedness  of idiosyncratic returns.

We focus on the S\&P 500 Financials sector. There are in total $80$ stocks on the first trading day of 2012 (Jan 3, 2012), among which $76$ stocks have complete data over the years
of 2012-2016. We will focus on these 76 stocks. The stock prices that our analysis is based on are collected from the Center for Research in Security Prices (CRSP) daily database,
 while the Fama-French three-factor data are obtained from Kenneth French's data library (\verb"http://mba.tuck.dartmouth.edu/pages/faculty/ken.french/data_library.html").

\begin{figure}[H]
	\begin{center}
		\includegraphics[width=3in]{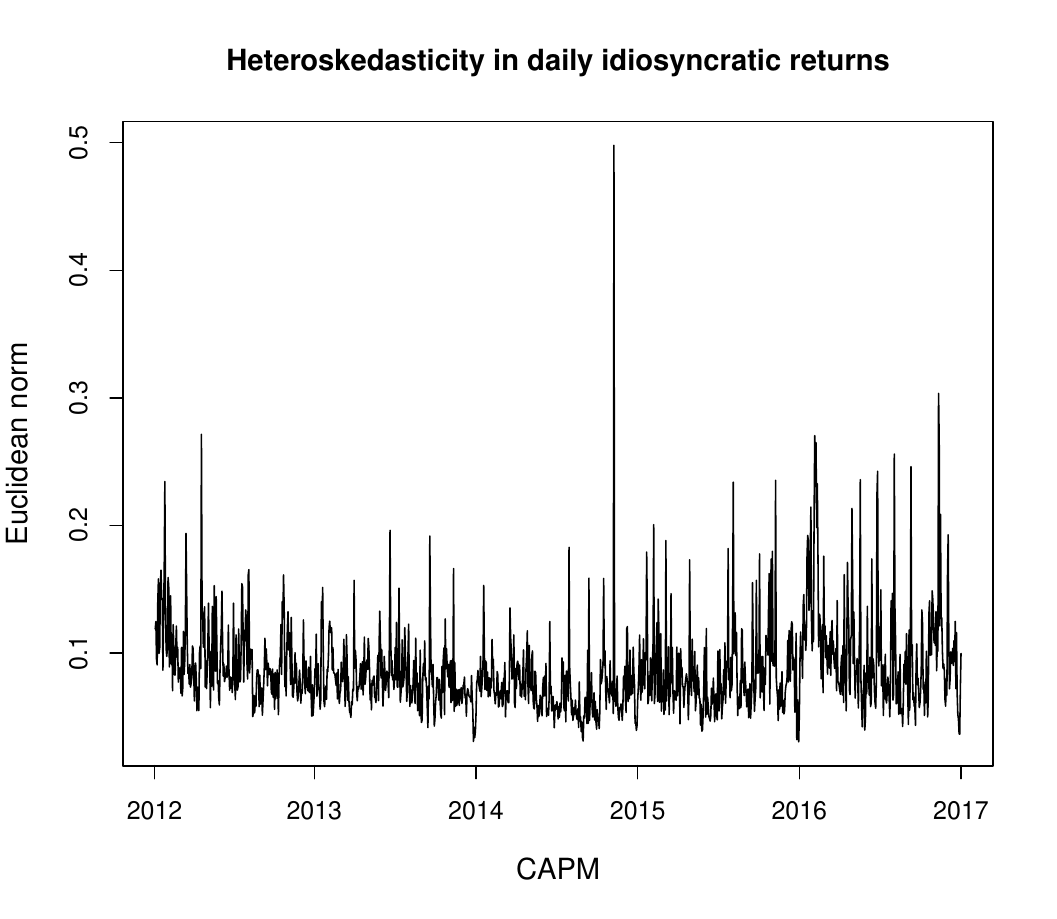}
        \includegraphics[width=3in]{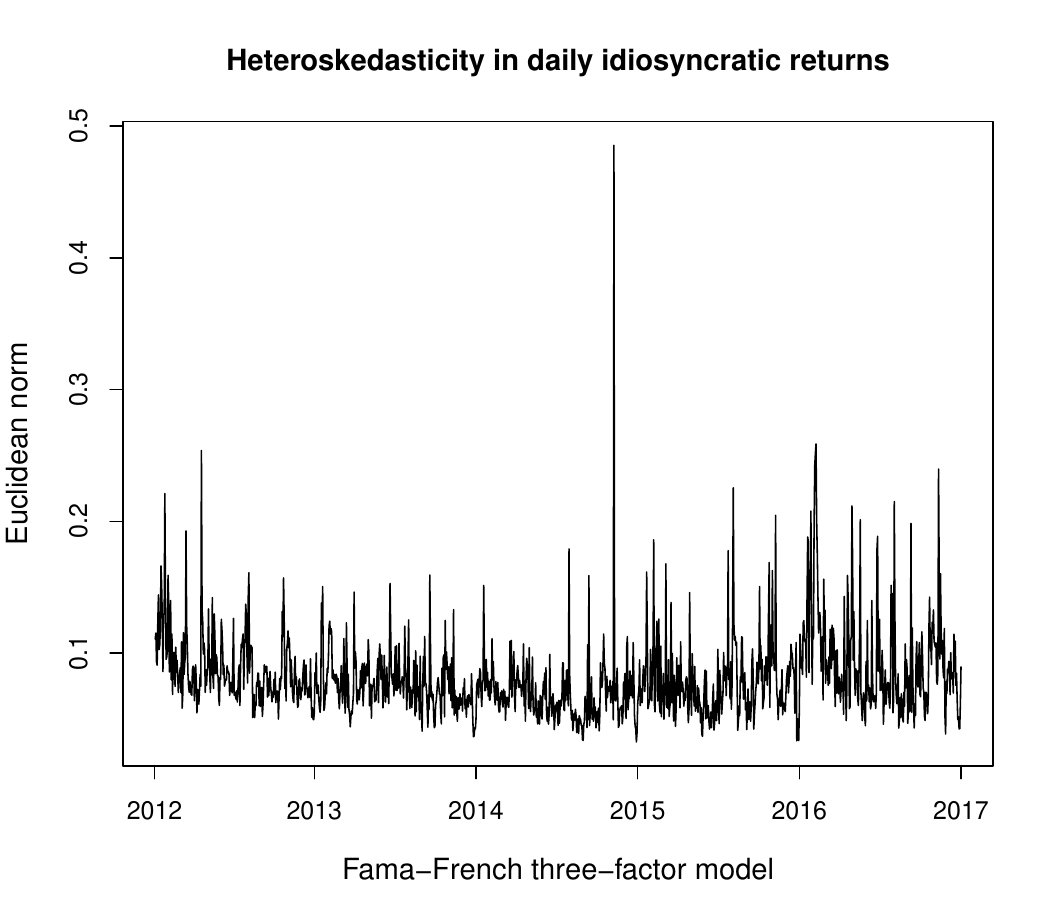}
	\end{center}
	\caption{}\label{fig:hetero}
\end{figure}

We illustrate the testing procedure based on two most widely used factor models: the CAPM and the Fama-French Three-Factor model.
We use a rolling window of six months to fit the two models.
Figure \ref{fig:hetero} reports the Euclidean norms of the fitted daily idiosyncratic returns.

\begin{center}
+++  Insert Figure \ref{fig:hetero} Here +++
\end{center}

We see from Figure \ref{fig:hetero} that under both models, the Euclidean norms of the fitted daily idiosyncratic returns exhibit clear heteroskedasticity and clustering. Such features indicate that the idiosyncratic returns are unlikely to be homoskedastic, but more suitably modeled as a conditional heteroskedastic time series, which is compatible with our framework.

Now we test the
 cross-sectional uncorrelatedness of idiosyncratic returns.
Specifically,
 for a diagonal matrix $\boldsymbol{\Sigma}_{\mathcal{D}}$ to be chosen, we test
\begin{align}
	&H_{0}:~\boldsymbol{\Sigma}_{\mathcal{I}} \propto \boldsymbol{\Sigma}_{\mathcal{D}}\quad vs. \quad H_{a}:~\boldsymbol{\Sigma}_{\mathcal{I}} \not\propto\boldsymbol{\Sigma}_{\mathcal{D}},\label{test:diag}
\end{align}
where $\boldsymbol{\Sigma}_{\mathcal{I}}$ denotes the covariance matrix of the  idiosyncratic returns. We will test \eqref{test:diag} by applying JHN-SN test to transformed idiosyncratic returns  by multiplying $\boldsymbol{\Sigma}_{\mathcal{D}}^{-1/2}$ to the (estimated) idiosyncratic returns.

\subsection{Testing results}\label{ssec:empirical_realdata}

We test \eqref{test:diag} using the same rolling window scheme as for fitting the CAPM or the Fama-French three-factor model. For each month to be tested, the
 diagonal matrix $\boldsymbol{\Sigma}_{\mathcal{D}}$ in \eqref{test:diag} is obtained by extracting the diagonal entries of the
sample covariance matrix of the self-normalized fitted idiosyncratic returns over the previous five months.
Table \ref{tab:t_diag} summarizes the resulting JHN-SN test statistics.

\begin{center}
+++  Insert Table \ref{tab:t_diag} Here +++
\end{center}

We observe from Table \ref{tab:t_diag} that:
\begin{enumerate}[(i)]
	\item The values of the JHN-SN test statistics are in general rather big, which correspond to {almost zero} $p$-values. Such a finding casts doubt on the cross-sectional uncorrelatedness of the idiosyncratic returns from fitting either the CAPM or the Fama-French three-factor model;
	\item Compared with the CAPM, the Fama-French three-factor model gives rise to idiosyncratic returns that are associated with less extreme test statistics. This confirms that the two additional factors, size and value, do have pervasive impacts on stock returns.
\end{enumerate}

\subsection{Robustness check of the testing results}

{The results in Table \ref{tab:t_diag} are based on testing against the estimated diagonal matrix $\boldsymbol{\Sigma}_{\mathcal{D}}$, which inevitably contains estimation errors. This brings up the following question: are the extreme test statistics in Table \ref{tab:t_diag} due to the estimation error in $\boldsymbol{\Sigma}_{\mathcal{D}}$, or, are they really due to that the idiosyncratic returns are not uncorrelated?}
To answer this question, we redo the test based on simulated stock returns whose idiosyncratic returns are uncorrelated and exhibit heteroskedasticity.

Specifically, we consider the following three-factor model:
 \begin{align}\label{three_factor_model}
\fr_t=\boldsymbol{\alpha}+\fB\ff_t+\boldsymbol{\varepsilon}_t,~~\text{with}~~ \ff_t\sim N(\boldsymbol{\mu}_f, \boldsymbol{\Sigma}_f),~~ \boldsymbol{\varepsilon}_t=\omega_t\cdot\boldsymbol{\Sigma}_{\boldsymbol{\mathcal{I}}}^{1/2}\fZ_t~~\text{and}~~\fZ_t\sim N(\boldsymbol{0}, \fI),
 \end{align}
 where $\fr_t$ denotes return vector at time $t$,
 $\fB$ is a factor loading matrix, $\ff_t$ represents three factors, and $\boldsymbol{\varepsilon}_t$ consists of idiosyncratic returns.
 To mimic the real data, we calibrate the parameters as follows:
  \begin{enumerate}[(i)]
  	\item The factor loading matrix  $\boldsymbol{B}$ is taken to be the estimated factor loading matrix by fitting the Fama-French three-factor model to the daily returns of the $76$ stocks over the years of 2012--2016, and $\boldsymbol{\alpha}$ is obtained by hard thresholding the estimated intercepts by two standard errors;
  	\item The mean and covariance matrix of factor returns, $\boldsymbol{\mu}_f$ and ${\fSigma}_f$, are the sample mean and sample covariance matrix of the Fama-French three factor returns from 2012 to 2016;
  	\item To generate data under the null hypothesis that the idiosyncratic returns are uncorrelated, {their covariance matrix $\boldsymbol{\Sigma}_{\boldsymbol{\mathcal{I}}}$
  is taken to be the diagonal matrix obtained by extracting the diagonal entries of the sample covariance matrix of the self-normalized fitted idiosyncratic returns}; and
  \item Finally, $\omega_t$ is taken to be the Euclidean norm of the fitted daily idiosyncratic returns.
  \end{enumerate}

With such generated data, we  test \eqref{test:diag} in parallel with the real data analysis.
Table \ref{tab:simu_diagonal} summarizes the JHN-SN test statistics for testing \eqref{test:diag} based on the simulated data.

\begin{center}
+++  Insert Table \ref{tab:simu_diagonal} Here +++
\end{center}

Table \ref{tab:simu_diagonal} reveals sharp contrast with Table \ref{tab:t_diag}. We see that if the idiosyncratic returns are indeed uncorrelated,
then even if they are heteroskedastic and even if we are testing against the estimated $\widehat{\boldsymbol{\Sigma}}_{\mathcal{D}}$, the percentage of resulting test statistics that are
within $[-1.96, 1.96]$ is close to 95\%, the expected level under the null hypothesis. In sharp contrast, the test statistics in Table \ref{tab:t_diag} are all very extreme.
 Such a comparison suggests that the idiosyncratic returns in the real data are indeed unlikely to be uncorrelated.

\subsection{Broader usage of the proposed test}

The testing procedure above can be directly translated to other factor models. Furthermore, the comparison made  in Section \ref{ssec:empirical_realdata} based on our test statistic can be viewed as a ``scoring'' system for different factor models. More specifically,  less extreme test statistic values would suggest that the model is more effective in accounting for pervasive impact of underlying factors.

\section{Conclusions}\label{conclusion}
We study  testing high-dimensional covariance matrices under a generalized elliptical distribution, which can feature heteroskedasticity, leverage effect, asymmetry, etc. We establish a CLT for the LSS of the sample covariance matrix based on self-normalized observations. The CLT is different from the existing ones for the usual sample covariance matrix. When the covariance matrix equals the identity matrix, our CLT neither requires $\mathbb{E}(Z_{11}^4)=3$ as in \cite{BaiS04} nor involves $\mathbb{E}(Z_{11}^4)$ as in \cite{Najim2013}. Based on the new CLT, we propose two sphericity tests by modifying the likelihood ratio test and John's test. {More general tests are also provided.} Numerical studies show that our proposed tests work well no matter whether the observations are \hbox{i.i.d.} Gaussian or from an elliptical distribution or feature conditional heteroskedasticity or  even when $\mathbf Z_i$'s do not admit the fourth moment. Moreover, we can also test against a general non-negative definite matrix which can be even not invertible.
As an innovative application, we demonstrate that our tests can be utilized to test uncorrelatedness among idiosyncratic returns. The testing procedure is illustrated by using CAPM and Fama-French Three-Factor model. The analysis can be translated directly to other factor models.


\newpage

\begin{table}[H]
		\centering
		\ra{0.63}\setlength{\tabcolsep}{3pt}	
		\begin{tabular}{@{}c||ccccccccccc@{}}
\toprule[1pt]
\multicolumn{11}{c}{$H_0: \boldsymbol{\Sigma}=\fI$}\\
\midrule[1pt]
			&\multicolumn{5}{c}{$p/n=0.5$} & \phantom{abc}& \multicolumn{4}{c}{$p/n=2$}\\
			\cmidrule{2-6} \cmidrule{8-11}
			$p$& LW$_1$ &BJYZ& CZZ$_1$ & WYMC-LR& WYMC-LW && LW$_1$ & CZZ$_1$&\multicolumn{2}{c}{WYMC-LW}\\ \midrule[1pt]
			$100$&$100$&$100$&$54.0$&$100$&$100$&&$100$&$50.2$&\multicolumn{2}{c}{$100$}\\
			$200$&$100$&$100$&$51.6$&$100$&$100$&&$100$&$53.0$&\multicolumn{2}{c}{$100$}\\
			$500$&$100$&$100$&$52.3$&$100$&$100$&&$100$&$53.3$&\multicolumn{2}{c}{$100$}\\
\midrule[1pt]
\multicolumn{11}{c}{$H_0: \boldsymbol{\Sigma}\propto\fI$}\\
\midrule[1pt]
			~~~~~~&\multicolumn{5}{c}{$p/n=0.5$} & \phantom{abc}& \multicolumn{4}{c}{$p/n=2$}\\
			\cmidrule{2-6} \cmidrule{8-11}
			$p$& LW$_2$ &S &  CZZ$_2$& WY-LR& WY-JHN && LW$_2$ & S& CZZ$_2$& WY-JHN \\ \midrule[1pt]
			$100$&$100$& $100$ & $51.8$&$100$& $100$ && $100$&$100$&$50.2$& $100$\\
			$200$&$100$& $100$& $53.0$& $100$& $100$ && $100$&$100$ &$52.3$& $100$\\
			$500$&$100$& $100$& $52.3$&$100$& $100$ && $100$&$100$&$53.5$& $100$\\
			\bottomrule[1pt]
		\end{tabular}
	\caption{\it Empirical sizes $(\%)$ of the existing tests for testing $H_0: \boldsymbol{\Sigma}=\fI$ or $H_0: \boldsymbol{\Sigma}\propto\fI$ at $5\%$ significance level.  Data are generated as $\fY_i=\omega_i\fZ_i$ where $\omega_i$'s are  absolute values of \hbox{i.i.d.} $N(0, 1)$, $\fZ_i$'s are \hbox{i.i.d.} $N(\mathbf{0}, \fI)$, and further $\omega_i$'s and $\fZ_i$'s are independent of each other. The results are based on $10,000$ replications for each pair of $p$ and $n$. }\label{counterexample_sphere}
\end{table}

\begin{table}[H]
\begin{center}
	\centering
	\ra{0.63}
	\setlength{\tabcolsep}{3.5pt}
	\begin{tabular}{@{}c||ccccc|ccccccc|c@{}}\toprule[1pt]
		~~~~~~&\multicolumn{7}{c}{$p/n=0.5$} & \phantom{abc}& \multicolumn{5}{c}{$p/n=2$}\\
		\cmidrule{2-8} \cmidrule{10-14}
		$p$& LW$_2$ &S &CZZ$_2$& WY-LR&LY& LR-SN & JHN-SN && LW$_2$ & S& CZZ$_2$& LY& JHN-SN \\ \midrule[1pt]
		$100$& $4.9$ & $4.8$&$4.9$ &$4.5$&$4.8$ & $4.6$& $5.2$&& $5.5$& $5.5$& $5.7$&$5.1$& $4.9$\\
		$200$&$5.2$& $5.0$& $5.1$& $5.1$&$4.8$& $5.1$& $4.9$&& $4.6$&$4.5$&$5.1$& $4.8$&$4.5$\\
		$500$&$4.9$& $5.1$&$5.1$&$4.8$&$5.3$&$4.9$&$5.2$ &&$5.1$&$5.3$&$4.9$&$5.0$&$5.2$\\
		\bottomrule[1pt]
	\end{tabular}
\end{center}
\caption{\it Empirical sizes $(\%)$ of LW$_2$, S, CZZ$_2$, WY-LR, LY, and the LR-SN and JHN-SN tests for testing $H_0: \boldsymbol{\Sigma}\propto\fI$ at $5\%$ significance level.  Observations are \hbox{i.i.d.}~$N(\mathbf{0}, \fI)$.  The results are based on $10,000$ replications for each pair of $p$ and $n$.}\label{compare_size}
\end{table}

\begin{table}[H]
\begin{center}
	\centering
	\ra{0.63}
	\setlength{\tabcolsep}{2.5pt}
	\begin{tabular}{@{}c||ccccc|ccccccc|c@{}}\toprule[1pt]
		~~~~~~&\multicolumn{7}{c}{$p/n=0.5$} & \phantom{abc}& \multicolumn{5}{c}{$p/n=2$}\\
		\cmidrule{2-8} \cmidrule{10-14}
		$p$& LW$_2$ &S &CZZ$_2$& WY-LR& LY& LR-SN & JHN-SN && LW$_2$ & S& CZZ$_2$ &LY& JHN-SN \\ \midrule[1pt]
		$100$& $50.7$ & $51.3$&$50.1$ &$36.7$&$28.1$ & $35.0$& $48.9$&& $8.4$& $8.7$& $9.1$&$6.3$& $8.2$\\
		$200$&$97.3$& $97.3$& $97.2$& $88.0$&$79.4$& $88.7$& $97.0$&&$18.3$& $17.9$&$18.1$&$11.9$&$17.2$\\
		$500$&$100$& $100$& $100$& $100$&$100$& $100$& $100$&& $70.7$&$70.6$&$69.8$&$43.3$& $70.5$\\
		\bottomrule[1pt]
	\end{tabular}
\end{center}
\caption{\it Empirical powers $(\%)$ of LW$_2$, S, CZZ$_2$, WY-LR, LY, and the LR-SN and JHN-SN tests for testing $H_0: \boldsymbol{\Sigma}\propto\fI$ at $5\%$ significance level. Observations are \hbox{i.i.d.} $N\Big(\mathbf{0}, \big(0.1^{|i-j|}\big)\Big)$. The results are based on $10,000$ replications for each pair of $p$ and $n$. }\label{compare_power}
\end{table}

\begin{table}[H]
	\begin{center}
		\ra{0.63}
		\setlength{\tabcolsep}{1.4pt}
		\begin{tabular}{@{}c||cccccc|cccccccc|c@{}}\toprule[1pt]
			~~~~~~&\multicolumn{8}{c}{$p/n=0.5$} & \phantom{~}& \multicolumn{6}{c}{$p/n=2$}\\
			\cmidrule{2-9} \cmidrule{11-16}
			$p$&{\small LW$\!_2$} &S &{\small CZZ$_2$}& {\small WY-LR}& {\small WY-JHN} &{\small LY}&{\small LR-SN}&{\small JHN-SN}&&{\small LW$\!_2$ }&{\small S}&{\small CZZ$_2$}& {\small WY-JHN}&{\small LY}&{\small JHN-SN}\\ \midrule[1pt]
			 {$100$}&{$100$}& $100$ & $51.8$&$100$& $100$&${4.4}$ &$\bf{4.6}$&$\bf{5.2}$ && $100$&$100$&$50.2$& $100$&${4.1}$&$\bf{4.9}$\\
			$200$&$100$& $100$& $53.0$& $100$& $100$ &${4.5}$&$\bf{5.1}$&$\bf{4.9}$ && $100$&$100$ &$52.3$& $100$&${4.5}$&$\bf{4.5}$\\
			$500$&$100$& $100$& $52.3$&$100$& $100$&${5.2}$&$\bf{4.9}$&$\bf{5.2}$ && $100$&$100$&$53.5$& $100$&${4.7}$&$\bf{5.2}$\\
			\bottomrule[1pt]
		\end{tabular}
	\end{center}
	\caption{\it{Empirical sizes $(\%)$ of LW$_2$, S, CZZ$_2$, WY-LR, WY-JHN, LY tests, and our proposed LR-SN, JHN-SN tests for testing $H_0: \boldsymbol{\Sigma}\propto\fI$ at $5\%$ significance level. Data are generated as $\fY_i=\omega_i\fZ_i$ where $\omega_i$'s are absolute values of \hbox{i.i.d.} $N(0, 1)$, $\fZ_i$'s are \hbox{i.i.d.}~$N(\mathbf{0}, \fI)$, and further $\omega_i$'s and $\fZ_i$'s are independent of each other. The results are based on $10,000$ replications for each pair of $p$ and~$n$. }}\label{ellipticalsize}
\end{table}

\begin{table}[H]
	\begin{center}
		\centering
		\ra{0.63}
		\begin{tabular}{@{}c||c|cccc|ccc|c@{}}\toprule[1pt]
			~~~~~~&\multicolumn{3}{c}{$p/n=0.5$} &\phantom{abc}& \multicolumn{2}{c}{$p/n=2$}\\
			\cmidrule{2-4} \cmidrule{6-7}
			$p$&  LY&LR-SN& JHN-SN&&LY& JHN-SN\\ \midrule[1pt]
			$100$&$7.6$&$35.0$& $48.9$&&$3.5$& $8.2$\\
			$200$&$14.5$&$88.7$& $97.0$&&$5.7$& $17.2$\\
			$500$&$64.9$&$100$& $100$&&$9.0$& $70.5$\\
			\bottomrule[1pt]
		\end{tabular}
	\end{center}
	\caption{\it{ Empirical powers $(\%)$ of LY test and our proposed LR-SN and JHN-SN tests   for testing $H_0: \boldsymbol{\Sigma}\propto\fI$ at $5\%$ significance level.
			Data are generated as $\fY_i=\omega_i\fZ_i$ where $\omega_i$'s are  absolute values of \hbox{i.i.d.} $N(0, 1)$, $\fZ_i$'s are \hbox{i.i.d.} random vectors from $N(\mathbf{0},  \boldsymbol{\Sigma})$ with $\boldsymbol{\Sigma}=\big(0.1^{|i-j|}\big)$, and further $\omega_i$'s and $\fZ_i$'s are independent of each other. The results are based on $10,000$ replications for each pair of $p$ and~$n$. }}\label{ellipticalpower}
\end{table}

\begin{table}[H]
	\begin{center}
		\centering
		\ra{0.63}
		\begin{tabular}{@{}c||c|cccc|ccc|c@{}}\toprule[1pt]
			~~~~~~&\multicolumn{3}{c}{$p/n=0.5$} & \phantom{abc}& \multicolumn{2}{c}{$p/n=2$}\\
			\cmidrule{2-4} \cmidrule{6-7}
			$p$&  LY&LR-SN& JHN-SN&&  LY& JHN-SN\\ \midrule[1pt]
			$100$&$8.2$&$\bf{5.5}$& $\bf{5.3}$&&$6.8$& $\bf{5.0}$\\
			$200$&$8.5$&$\bf{5.7}$& $\bf{5.4}$&&$6.8$& $\bf{5.5}$\\
			$500$&$7.6$&$\bf{5.3}$& $\bf{5.2}$&&$6.6$& $\bf{5.4}$\\
			\bottomrule[1pt]
		\end{tabular}
	\end{center}
	\caption{\it Empirical sizes $(\%)$ of LY test and our proposed LR-SN and JHN-SN tests   for testing $H_0: \boldsymbol{\Sigma}\propto\fI$ at $5\%$ significance level. Data are generated as $\fY_i=\omega_i\fZ_i$ with $\omega^2_i=0.01+0.85\omega^2_{i-1}+0.1|\fY_{i-1}|^2/p$, and $\fZ_i$'s consist of \mbox{i.i.d.} standardized $t(4)$ random variables. The results are based on $10,000$ replications for each pair of $p$ and $n$. }\label{size}
\end{table}

\begin{table}[H]
	\begin{center}
		\centering
		\ra{0.63}
		\begin{tabular}{@{}c||c|cccc|ccc|c@{}}\toprule[1pt]
			~~~~~~&\multicolumn{3}{c}{$p/n=0.5$} &  \phantom{abc}& \multicolumn{2}{c}{$p/n=2$}\\
			\cmidrule{2-4} \cmidrule{6-7}
			$p$&  LY&LR-SN& JHN-SN&&  LY& JHN-SN\\ \midrule[1pt]
			$100$&$20.7$&$34.4$& $47.9$&&$7.8$& $8.7$\\
			$200$&$54.4$&$87.8$& $96.6$&&$10.5$& $17.6$\\
			$500$&$100$&$100$& $100$&&$26.4$& $69.9$\\
			\bottomrule[1pt]
		\end{tabular}
	\end{center}
	\caption{\it Empirical powers $(\%)$ of LY test and our proposed LR-SN and JHN-SN tests  for testing $H_0: \boldsymbol{\Sigma}\propto\fI$ at $5\%$ significance level.
		Data are generated as $\fY_i=\omega_i\fSigma^{1/2}\fZ_i$ with $\omega^2_i=0.01+0.85\omega^2_{i-1}+0.1|\fY_{i-1}|^2/p$ and $\fSigma=\big(0.1^{|i-j|}\big)$, and ${\fZ}_i$’s consist of \mbox{i.i.d.} standardized $t(4)$ random variables. The results are based on $10,000$ replications for each pair of $p$ and~$n$. }\label{tab:power}
\end{table}

\begin{table}[H]
	\begin{center}
		\centering
		\ra{0.63}
		\begin{tabular}{@{}c||ccc@{}}\toprule[1pt]
			$p$ &\multicolumn{1}{c}{$p/n=0.5$} & \phantom{abc}& \multicolumn{1}{c}{$p/n=2$}\\
			 \midrule[1pt]
			$100$& $4.8$&&$4.4$\\
			$200$& $4.9$&&$4.7$\\
			$500$& $5.1$&&$4.6$\\
			\bottomrule[1pt]
		\end{tabular}
	\end{center}
	\caption{\it Empirical sizes $(\%)$ of our test for testing $H_0: \boldsymbol{\Sigma}\propto\boldsymbol{\Sigma}_0$ at $5\%$ significance level. Here  $\boldsymbol{\Sigma}_0=\mathbf{Q}\boldsymbol{\Lambda}_0\mathbf{Q}^\top$, where $F^{\boldsymbol{\Lambda}_0}=\frac{1}{2}\delta_1+\frac{1}{4}\delta_2+\frac{1}{4}\delta_0$ 	
and $\fQ$ is a random orthogonal matrix. Data are generated as $\fY_i=\omega_i\fSigma_0^{1/2}\fZ_i$ with $\omega^2_i=0.01+0.85\omega^2_{i-1}+0.1|\fY_{i-1}|^2/p$ and $\fZ_i$'s from $N(\mathbf{0}, \fI)$. The results are based on $10,000$ replications for each pair of $p$ and $n$.
}\label{tab:size_general}
\end{table}

\begin{table}[H]
	\begin{center}
		\centering
		\ra{0.63}
		\begin{tabular}{@{}c||ccc@{}}\toprule[1pt]
			$p$ &\multicolumn{1}{c}{$p/n=0.5$} & \phantom{abc}& \multicolumn{1}{c}{$p/n=2$}\\
			 \midrule[1pt]
			$100$& $100$&&$46.3$\\
			$200$& $100$&&$95.5$\\
			$500$& $100$&&$100$\\
			\bottomrule[1pt]
		\end{tabular}
	\end{center}
	\caption{\it Empirical power $(\%)$ of our test for testing $H_0: \boldsymbol{\Sigma}\propto\fSigma_0$ at $5\%$ significance level. Data are generated as $\fY_i=\omega_i\fSigma^{1/2}\fZ_i$ with $\omega^2_i=0.01+0.85\omega^2_{i-1}+0.1|\fY_{i-1}|^2/p$ and $\fZ_i$'s from $N(\mathbf{0}, \fI)$. Here  $\fSigma=\fQ\boldsymbol{\Lambda}\fQ^\top$, where $F^{\boldsymbol{\Lambda}}=\frac{1}{2}\delta_1+\frac{1}{4}\delta_{2.5}+\frac{1}{4}\delta_0$   and $\fQ$ is a random orthogonal matrix.
 The results are based on $10,000$ replications for each pair of $p$ and $n$. }\label{tab:power_general}
\end{table}

\begin{table}[H]
	\begin{center}
		\ra{0.6}
		\begin{tabular}{cccccc|c}
			\toprule[1pt]
			\multicolumn{7}{c}{CAPM}\\
			\toprule[1pt]
			~&Min& $Q_1$& Median& $Q_3$& Max& Mean (Sd)\\
			JHN-SN& $6.3$& $18.1$& $29.8$& $44.3$& $83.1$&$33.1~(18.5)$\\
			\toprule[1pt]
			\multicolumn{7}{c}{Fama-French three-factor model}\\
			\toprule[1pt]
			~&Min& $Q_1$& Median& $Q_3$& Max& Mean (Sd)\\
			JHN-SN& $5.0$& $12.4$& $24.4$& $30.4$& $77.0$&$23.8~(13.0)$\\		
			\toprule[1pt]
		\end{tabular}
	\end{center}
\caption{\it Summary statistics of the JHN-SN statistics  for testing \eqref{test:diag}. For both the CAPM and the  Fama-French three-factor model, for each month,
we first estimate the idiosyncratic returns by fitting the model using the data in the current month and the previous five months. We then obtain
 $\boldsymbol{\Sigma}_{\mathcal{D}}$ by extracting the diagonal entries of the sample covariance matrix of
the self-normalized idiosyncratic returns over the previous five months, and use the fitted idiosyncratic returns in the current month to conduct the test.}\label{tab:t_diag}
\end{table}

\begin{table}[H]
 \ra{0.7}
	\begin{center}
		\begin{tabular}{cccccc|c|c}
			\toprule[1pt]
			\multicolumn{8}{c}{Simulated data based on a three-factor model}\\
			\toprule[1pt]
			~&Min& $Q_1$& Median& $Q_3$& Max& Mean (Sd)& Percent within $[-1.96, 1.96]$ \\
			JHN-SN& $-1.1$& $-0.2$& $0.6$& $1.2$& $2.5$&$0.6~(0.9)$& $94.5\%$\\				
			\bottomrule[1pt]
		\end{tabular}
	\end{center}
	\caption{\it Summary statistics of the JHN-SN statistics  for testing \eqref{test:diag} based on simulated returns from Model \eqref{three_factor_model}. {To conduct the test, with a rolling window of six months,
 we first estimate the idiosyncratic returns by fitting the three-factor model.} We then obtain
 $\boldsymbol{\Sigma}_{\mathcal{D}}$ by extracting the diagonal entries of the sample covariance matrix of
the self-normalized fitted idiosyncratic returns over the previous five months, and use the fitted idiosyncratic returns in the current month to conduct the test.}\label{tab:simu_diagonal}
\end{table}

\newpage
\section*{Figure legends}
\noindent {\bf Fig.1} Time series plots of the Euclidean norms of the daily idiosyncratic returns of $76$ stocks in the S$\&$P 500 Financials sector, by fitting the CAPM (left) and the Fama-French three-factor model (right) over the years of 2012--2016.

\newpage

\end{document}